\documentclass[12pt]{article}
\newcommand{\floor}[1]{\left\lfloor #1 \right\rfloor}
\newcommand{\ffloor}[2]{\left\lfloor \frac{#1}{#2} \right\rfloor}

\renewcommand{\r}[2][p]{r_{#1}\left(#2\right)}
\newcommand{\rk}{r_{p}\left(k^2\right)}
\newcommand{\rkm}[1][k]{r_{p}\left((#1-1)^2\right)}
\newcommand{\rkl}{r_{p}\left(k+1-3n\right)}

\newcommand{\kt}[1][m]{(k-1)^2-#1p+k+1-3n}

\newcommand{\keh}[1][m]{(k-1)^2-#1p+n-k+1-p}
\newcommand{\khhh}[1][m]{(\widehat{k}-1)^2-\widehat{#1}p+n-\widehat{k}+1-p}
\newcommand{\khkh}[1][m]{(k-1)^2-\widehat{#1}p+n-\widehat{k}+1-p}

\newcommand{\fq}[1]{\left\lfloor Q_{#1} \right\rfloor}
\newcommand{\fqm}{\left\lfloor Q_m \right\rfloor}

\newcommand{\fr}[1]{\left\lfloor R_{#1} \right\rfloor}
\newcommand{\dfr}[1]{\left\lfloor \newsqrt{#1p} \right\rfloor}
\newcommand{\dfrk}{\left\lfloor \newsqrt{kp} \right\rfloor}
\newcommand{\frm}{\left\lfloor R_m \right\rfloor}

\newcommand{\fs}[1]{\left\lfloor S_{#1} \right\rfloor}

\newcommand{\rp}[1][1]{\frac{#1+\sqrt{4mp+12n-7}}{2}}
\newcommand{\rpi}[1][1]{\left(#1+\sqrt{4mp+12n-7}\right)/2}
\newcommand{\trpi}[1][1]{#1+\sqrt{4mp+12n-7}}
\newcommand{\yp}[2][1]{\frac{#1+\sqrt{4#2p+12n-7}}{2}}

\newcommand{\den}{4mp+12n-7}
\newcommand{\gk}{\Gamma(k)}
\newcommand{\gm}[1][n+1]{\Gamma(#1)}
\newcommand{\gpk}{\Gamma(p+2-k)}

\newcommand{\hAge}{C_{\ge}}
\newcommand{\hAbet}{C_{[\hspace{-.075cm}-\!\!-\!\!)}}
\newcommand{\hAle}{C_{<}}
\newcommand{\hBge}{D_{\ge}}
\newcommand{\hBbet}{D_{[\hspace{-.075cm}-\!\!-\!\!)}}
\newcommand{\hBle}{D_{<}}
\newcommand{\mm}[1]{$\,\ds{#1}\,$}
\newcommand{\m}[1]{$\,#1\,$}

\newcommand{\wk}{{\widehat{k}}}
\newcommand{\wh}{{\widehat{m}}}
\newcommand{\ns}{\hspace{-1.8cm}}

\newlength{\largo}
\newlength{\corto}
\newlength{\resta}

\newcommand{\quadf}{\Q(\sqrt{-p})}
\usepackage{float}
\restylefloat{table}

\newcommand {\ds}{\displaystyle}

\newtheorem{thm}{Theorem}[section]
\newtheorem{lem}[thm]{Lemma}
\newtheorem{obs}[thm]{Observation}
\newtheorem{cor}[thm]{Corollary}

\newtheorem{rem}[thm]{Remark}
\newtheorem{nota}[thm]{Notation}

\newtheorem{cnjt}[thm]{Conjecture}

\newtheorem{dfn}[thm]{Definition}

\newcommand{\ok}{\hfill $\Box$\\[2mm]}

\newlength{\tamano}
\newlength{\taman}
\newcounter{figura}
\newcommand{\mytext}[1]{\mbox{\ #1\ }}

\newcommand{\crea}[2]{\parbox[l][#2cm][l]{#1cm}{\ }}
\newcommand{\at}[2][rrrrrrrrrrrrrr]{    %Se debe utilizar así \at[llll]{1&2&3\\ 4&5&5}
                     \begin{array}{#1}
                     #2\\
                     \end{array}}
\newcommand{\defi}[2]{#1 = \left\{\at[lllllllllllll]{#2}\right.} %\defi{f(x)}{3 & x\le 4\\ 2 & x>4}

\def\newsqrt{\mathpalette\DHLhksqrt} % it defines the new \sqrt in terms of the old one
\def\DHLhksqrt#1#2{\setbox0=\hbox{$#1\sqrt{#2\,}$}\dimen0=\ht0
\advance\dimen0-0.45\ht0
\setbox2=\hbox{\hspace{-.014cm}\vrule height\ht0 depth -\dimen0}%
{\box0\lower0.55pt\box2}\,}

\newcommand{\nn}{\noindent}

\newcommand{\Q}{{\mathbb Q}}
\newcommand{\N}{{\mathbb N}}

\newcommand{\Z}{{\mathbb Z}}
\newcommand{\zz}{\Z_0}
\newcommand{\too}{\longrightarrow}
\renewcommand{\to}{\rightarrow}

\usepackage{amsfonts,amsmath,enumerate}

 \usepackage{graphicx}
 \usepackage{color}
 \usepackage{enumerate}
 \usepackage{epsfig}
 \usepackage{fancyhdr}
 \usepackage{ifthen}
 \usepackage{mathtools}
 \usepackage[tiling]{pst-fill}
 \usepackage{xifthen}
 \usepackage{xparse}
 \usepackage{tikz}
 \usepackage{float}

\usepackage{xcolor}
\usepackage{tabularray}
\begin{document}
\title{Class Number of the Imaginary Quadratic Field and Quadratic Residues Identities}
\author{Jorge Garcia}
\maketitle

%%%%%%%%%%%%%%%%%%%%%%%%%%%%%%%%%%%%%%%%%%%%%%%%%%%%%%%%%%%%%%%%%%%%%%%%%%%%%%%%%%%%%%%%%%%%%%%%%%%%%%%%%%%%%%%%%%%%%%%%%%%%%%%%%%%%%%%%%%%%%%%%%%%%%%%%%%%%%%%%%%%%%%%%%%%%%%%%%%%%%%%%%%%%%%%%%%%%%%%%%%%%%%%%%%%%%%%%%%%%%%%%%%%%%%%%%%%%%%%%%%%%%%%%%%%%%%%%%%%%%%%%%%%%%%%%%%%%%%%%%%%%%%%%%%%%%%%%%%%%%%%%%%%%%%%%%%%%%%%%%%%%%%

\begin{abstract}
   A formula for the sum of quadratic residues modulus a prime $p=4n-1$ is studied. We relate some terms on this formula with roots of quadratics and provide an exhaustive analysis of new concepts based on these roots. A number of formulas for the sum of the quadratic residues are obtained. We finalize the paper by obtaining several identities involving $h(-p)$ the class number of the imaginary quadratic field $\Q(\sqrt{-p}).$
\end{abstract}

%\begin{keywor d}
%  quadratic residues, prime numbers, sums, class number.
%\end{keyword}

%%%%%%%%%%%%%%%%%%%%%%%%%%%%%%%%%%%%%%%%%%%%%%%%%%%%%%%%%%%%%%%%%%
%%%%%%%%%%%%%%%%%%%%%%%%%%%%%%%%%%%%%%%%%%%%%%%%%%%%%%%%%%%%%%%%%%
%%%%%%%%%%%%%%%%   Section Introduction                    %%%%%%%
%%%%%%%%%%%%%%%%%%%%%%%%%%%%%%%%%%%%%%%%%%%%%%%%%%%%%%%%%%%%%%%%%%
%%%%%%%%%%%%%%%%%%%%%%%%%%%%%%%%%%%%%%%%%%%%%%%%%%%%%%%%%%%%%%%%%%

\section{Introduction}\label{Section:IntroAndDefinitions}
Consider a prime $p=4n-1$ and $1\le k\le p-1.$ By $\r{k^2},$ we denote the remainder of $k^2$ when we divide by $p$.  We call this number $\r{k^2}$ the \textit{quadratic residue} of $k^2$ modulus $p.$ When we add all these residues we obtain the sum of quadratic residues relative to the prime number $p$. There is a complicated formula for such sum,

 \begin{align}\label{DirichletFormula}
      \sum_{k=1}^{p-1}\r{k^2} =\binom{p}{2}-p\cdot h(-p),
    \end{align}

where $h(-p)$ is the class number of the imaginary quadratic field $\quadf$. Important formulas for the class number when the prime is of the form $p=4n-1>3$ include \textit{Dirichlet class number formula}

 \begin{eqnarray}\label{DirichletClassNumberFormula}
     % \nonumber to remove numbering (before each equation)
       h(-p) &=& \frac{\sqrt{p}}{2\pi}\sum_{r=1}^{\infty}\frac{\chi(r)}{r},
 \end{eqnarray}

where $\chi$ is the Dirichlet character.

Formula~\ref{DirichletClassNumberFormula} can be found in~\cite{Dirichlet:bedesat37}. Another formula that involves the Kronecker symbol $\Big(\frac{p}{r}\Big)$ can be found in \cite{Narkiewicz:eatan04} (Corollary 1 p. 428). Finally a formula developed by Cohen~\cite{cohen:ccant93} [Corollary 5.3.16] provides the class number too.

The main purpose of this paper is to obtain some identities for this number $h(-p)$ by computing the quadratic residues in a different manner.

Here we provide a summary of the main formulas obtained in this paper. Here $Q_k=\left(1+\newsqrt{4kp+12n-7}\right)/2,\ M=\floor{n^2/p}$\, and  $J_n$ is the number of jumps (see Section \ref{Section:SumOfQuadraticResidues}).

\begin{eqnarray*}
% \nonumber to remove numbering (before each equation)
    \sum_{k=1}^{M} \dfr{k}- \sum_{k=0}^{M-1}\fq{k} & = & J_n-M-1. \\
    \sum_{k=1}^{n}\ffloor{k^2-k+2-3n}{p}+\sum_{k=0}^{M-1}\floor{Q_k} & = & (M-1)n.
\end{eqnarray*}

\begin{eqnarray*}
    \sum_{k=1}^{p-1}\r{k^2-k+2-3n} & = & \sum_{k=1}^{p-1}\r{k^2}+3n-2.\\
    \sum_{k=1}^{2n}\r{k^2-k+2-3n}   & = & \sum_{k=1}^{2n}\r{k^2}+n.\\
    \sum_{k=1}^{n}\r{k^2-k+2-3n}   & = & \sum_{k=1}^{n}\r{k^2}+\frac{n(n+1)}{2}-p(J_n-1-M).\\
   \sum_{k=1}^{p-1}\ffloor{k^2-k+2-3n}{p} & = & \frac{p(p-5)+6-n}{3}-\frac{1}{p}\sum_{k=1}^{p-1}\rk.\\
   \sum_{k=1}^{2n}\ffloor{k^2-k+2-3n}{p} & = & \frac{n(2n-1)(4n-7)}{3p}-\frac{1}{p}\sum_{k=1}^{2n}\rk.\\
   \sum_{k=1}^{n}\ffloor{k^2-k+2-3n}{p} & = & \frac{(n+1)(2n^2-23n+6)}{3p}-\frac{1}{p}\sum_{k=1}^{n}\rk+J_n-M.
\end{eqnarray*}

This is a summary of the identities involving the class number $h=h(-p)$ found on this paper.

\begin{eqnarray*}
   \frac{J_n}{2}+M(n-1) & = &  \frac{h}{4}+\sum_{k=1}^{M}\dfrk+\frac{n^2-5n+9}{12}. \\
   -\frac{J_n}{2}+Mn & = &  \frac{h}{4}+\sum_{k=0}^{M-1}\fq{k}+\frac{n^2-5n-3}{12}.\\
   \frac{J_n}{2}+\sum_{k=1}^{n-1}\ffloor{k^2}{p} & = &  \frac{h}{4}+\frac{n^2-5n+9}{12}.
 \end{eqnarray*}

\begin{eqnarray*}
    \sum_{k=1}^{p-1}\r{k^2-k+2-3n} & = & p\cdot\left(2n-h\right)-n-1.\\
    \sum_{k=1}^{2n}\r{k^2-k+2-3n}   & = & \frac{p\cdot(2n-h)+1}{2}.\\
    \sum_{k=1}^{n}\r{k^2-k+2-3n}   & = & \frac{p}{4}(3n+2-2J_n-h)-\frac{(n+1)(n-1)}{4}.
\end{eqnarray*}

\begin{eqnarray*}
  \sum_{k=1}^{n}\r{k^2} & = & \frac{p}{4}(2J_n+2n-3-4M-h)+\frac{n(n+1)}{4}.
\end{eqnarray*}

\begin{eqnarray*}
   %\sum_{k=0}^{n}\ffloor{k^2-k+2-3n}{p} &=& \frac{h}{2}+\sum_{k=0}^{M}\fr{k}+M(1-n)+\frac{n^2-11n-3}{6},\\
   \sum_{k=1}^{p-1}\ffloor{k^2-k+2-3n}{p} &=& h+\frac{16n^2-35n+15}{3}.\\
   \sum_{k=1}^{2n}\ffloor{k^2-k+2-3n}{p} &=& \frac{h}{2}+\frac{4n^2-14n+9}{6}.\\
   \sum_{k=1}^{n}\ffloor{k^2-k+2-3n}{p} &=& \frac{h}{4}+\frac{J_n}{2}+\frac{n^2-17n-3}{12}.\\
  %\sum_{k=1}^{n}\ffloor{k^2-k+2-3n}{p} &=& \frac{h}{2}+M+\frac{n^2-11n+3}{6}- \sum_{k=1}^{n}\ffloor{k^2}{p}.\\
  \sum_{k=1}^{n}\ffloor{k^2-k+2-3n}{p} &=& \frac{h}{2}+M(1-n)+\frac{n^2-11n+3}{6}+\sum_{k=1}^{M}\dfrk.
 \end{eqnarray*}

\medskip

In Section \ref{Section:SumOfQuadraticResidues} we provide the main definitions and the notation used in the whole paper. Here, we develop the concepts of \textit{jump} and \textit{total residue}. We also identify when these jumps occur.  In Section \ref{Section:SplittingTheJumps} we organize the jumps on six different sets and we state the main lemmas that will be used to count the jumps which is done in Section \ref{Section:CountingTheJumps}. It is here where we define bijective functions among different pairs of sets to compute their cardinalities. In Section~\ref{Section:ResiduesInvolvingk^2-k+2-3n} we obtain the sums of quadratic residues of terms of the form $k^2-k+2-3n$ where $k$ ranges on the different intervals $[1,4n-2], [1,2n]$ or $[1,n].$  In this section we also count to total amount of jumps. Finally in Section \ref{Section:ClassNumber} we establish several identities involving the class number $h(-p)$ of the imaginary quadratic field $\Q(\sqrt{-p}).$ These identities are based on the sums found in previous sections and in some of these identities the jumps quantity appears.

Whereas in~\cite{garcia:siqrmp22} we computed several sums of quadratic residues when $p=4n+1$, in this paper, we perform a similar but different analysis when $p=4n-1.$ It is on this paper where the class number is involved. On some occasions we present both formulas ($p=4n-1$ and $p=4n+1$) for comparison purposes.

%%%%%%%%%%%%%%%%%%%%%%%%%%%%%%%%%%%%%%%%%%%%%%%%%%%%%%%%%%%%%%%%%%
%%%%%%%%%%%%%%%%%%%%%%%%%%%%%%%%%%%%%%%%%%%%%%%%%%%%%%%%%%%%%%%%%%
%%%%%%%%%%%%%%%%   Section Basic Lemmas                    %%%%%%%
%%%%%%%%%%%%%%%%%%%%%%%%%%%%%%%%%%%%%%%%%%%%%%%%%%%%%%%%%%%%%%%%%%
%%%%%%%%%%%%%%%%%%%%%%%%%%%%%%%%%%%%%%%%%%%%%%%%%%%%%%%%%%%%%%%%%%

\section{Sum of Quadratic Residues and Jumps}\label{Section:SumOfQuadraticResidues}

\nn Consider $p=4n-1$ a prime number. It will be understood that when we write $n$, we mean a natural number $n\ge 1.$

%%%%%%%%%%%%%%%%%%%%%%%%%%%%%%%%%%%%%%%%%%       Definition        %%%%%%%%%%%%%%%%
\begin{dfn}\label{remainder}
  Let $q$ be a positive integer and $x\in\Z$. By $r_q(x)$ we denote the remainder of $x$ when we divide by $q$. Hence $r_q(x)\in\{0,1,2,...,q-1\}$ satisfies
  \[x=m\cdot q + r_q(x),\]
  for some $m\in \Z$. Clearly, $m=\floor{x/q}.$
\end{dfn}

\nn The following notation found in~\cite{garcia:cfcniqf21} will be useful during the whole paper.

%%%%%%%%%%%%%%%%%%%%%%%%%%%%%%%%%%%%%%%%%%       Notation         %%%%%%%%%%%%%%%%
\begin{nota}\label{NotationQmRmM} For $m\in\Z, \,p=4n-1$ prime and $m\ge 0$ we denote
  \[
          \begin{split}
            Q_m & = \frac12+\frac12\newsqrt{1+4\left[(m+1)p-n-1\right]} = \frac12+\frac12\newsqrt{4mp+3p-4}, \\
                               R_m  & = \newsqrt{mp}\crea{0}{1.25},\\
                               M  & = \floor{\frac{n^2}{p}}.
          \end{split}
       \]

\end{nota}

In ~\cite{garcia:cfcniqf21}, we obtained a theorem that involves the sum of quadratic residues when $p=4n-1$.  For reference purposes, we write here such theorem.

%%%%%%%%%%%%%%%%%%%%%%%%%%%%%%%%%%%%%%%%%%       Corolario  9.5       %%%%%%%%%%%%%%%%
\begin{thm}\label{Corollary9.5TheoremMain4n-1}
Let $p=4n-1$ be prime. Using Notation~\ref{NotationQmRmM} we have
  \begin{equation*}
     \frac12\sum_{k=1}^{p-1}\r{k^2} = p\left(\sum_{m=1}^{M} \frm +\sum_{m=0}^{M-1} \fqm\right)-Mp(2n-1)+\frac{p\cdot(n^2+n)}{6}.
  \end{equation*}
\end{thm}

A concept arises naturally when we study the term $Q_m$. This term is the positive root of the quadratic polynomial $x^2-x+2-3n-mp$ which is the same as $(x-1)^2+x+1-3n-mp.$

%---------------------------------------------------------------------------------------------------------------
%%%%%%%%%%%%%%%%%%%%%%%%%%%%%%%%%%%%%%%%%%       Definition        %%%%%%%%%%%%%%%%
\begin{dfn}\label{totalresidue}
   Let $p=4n-1$ be a prime and $0\le k\le p$. The total residue of $k$ is defined and denoted by
   \[\Gamma(k)=\r{(k-1)^2}+\r{k+1-3n}.\]
   We also say that $k$ is a jump if its total residue is $p$ or more, i.e. if
   \[\Gamma(k)\ge p.\]
\end{dfn}
%----------------------------------------------------------------------------------------------------------------

What is the importance of the jumps? Firstly, in the proof of Theorem~\ref{Corollary9.5TheoremMain4n-1} (see~\cite{garcia:cfcniqf21}), a key technique is adding $\r{k^2}$ and $\r{(2n-k)^2}.$ It happens that when $k\in (R_m,Q_m]$ this sum is constant, but as soon as $k$ exceeds $Q_m$, we need to subtract $p.$ Therefore \textit{knowing when we need to subtract $p$} is key to comprehend better the formula in Theorem~\ref{Corollary9.5TheoremMain4n-1}. Secondly, the amount of jumps in the interval $[2,n+2]$ will allow us to establish a formula to compute the terms $\sum_{m=0}^{M-1}\floor{Q_{m}},\ \sum_{m=1}^{M}\floor{R_{m}}$ as well as $\sum_{k=0}^{n}\r{k^2}$ as a function of $n, h(-p)$ and the number of jumps. This is achieved in Corollaries~\ref{Corollary75} and~\ref{Corollary77}.

\medskip

The following three lemmas allow us to identify some jumps and when they occur.
For notation purposes we define $\zz=\{0,1,2,...\}.$

%%%%%%%%%%%%%%%%%%%%%%%%%%%%%%%%%%%%%%%%%%       Lemma 53             %%%%%%%%%%%%%%%%
\begin{lem}\label{Lemma53QhInequality}
  Let $p=4n-1$ be a prime and $m\in\zz.$  Then
 \[ \fqm<Q_m.\]
\end{lem}
%%%%%%%%%%%%%%%%%%%%%%%%%%%%%%%%%%%%%%%%%%       Proof of Lemma    %%%%%%%%%%%%%%%%
\noindent {\bf Proof.}
Assume there is $j\in\Z$ such that \mm{j=\rpi.} Then \m{j^2-j+2-3n=mp.} Taking \m{x_0=2n-j} gives
\[x_0^2=4n^2-n+mp+j-4nj+n+3n-2=(n+m-j+1)p-1,\]
hence \m{x_0^2\equiv -1 \pmod p} which contradicts Fermat Little Theorem as $p=4n-1.$ \ok

%%%%%%%%%%%%%%%%%%%%%%%%%%%%%%%%%%%%%%%%%%       Remark 53.5          %%%%%%%%%%%%%%%%
\begin{rem}\label{Note53.5k^2-3k+4-3n=hp}
By using the same argument, there is no integer $k$ with $k^2-3k+4-3n=mp,$ else taking $k-1=j$ we obtain \m{j^2-j+2-3n=mp} which leads to a contradiction.
\end{rem}

%%%%%%%%%%%%%%%%%%%%%%%%%%%%%%%%%%%%%%%%%%       Lemma  54           %%%%%%%%%%%%%%%%
\begin{lem}\label{Lemma54k=Qh-isajump}
  Let $p=4n-1$ be a prime with $n\ge 3.$ \ Let $m\in\zz$ with
   \[0\le m\le \ffloor{n^2-4n+5}{p} \textit{\ \ and\ \ } k_m=1+\floor{\rp}.\]
   Then
  \begin{enumerate}[(i)\,]
    \item \m{3\le k_m\le n+2.}
    \item \m{\trpi<2k_m<\trpi[3].}
    \item $k_m$ is a jump and
    \item \mm{3n-k_m-1<\r{(k_m-1)^2}=(k_m-1)^2-mp\le3n+k_m-4.}
  \end{enumerate}
\end{lem}

%%%%%%%%%%%%%%%%%%%%%%%%%%%%%%%%%%%%%%%%%%       Proof of Lemma    %%%%%%%%%%%%%%%%
\noindent {\bf Proof.}
\begin{enumerate}[(i)\,]
  \item Notice that \m{\den\le 4n^2+4n+1,} hence \m{k_m\le n+2.} Since $2\le n,\   9\le\den.$ Therefore  $3\le k_m.$
  \item Notice that $k_m$ is strictly above the positive root of $x^2-x+2-3n-mp,$ hence
      \begin{equation}\label{Lemma54Inequality1}
        k_m^2-k_m+2-3n>mp.
      \end{equation}

      Now \m{k_m\le n+2\le 3n} and \m{2\le n} imply
      \begin{equation}\label{Lemma54Inequality2}
        (k_m-1)^2-mp=k_m^2-k_m+2-3n-mp+3n-k_m-1\ge 3n-k_m.
      \end{equation}
      From Lemma~\ref{Lemma53QhInequality}, \m{k_m<\rpi[3].} The other inequality is obvious.
  \item
      From (ii), $k_m$ is less than the positive root of $x^2-3x+4-mp-3n$, hence  $k_m^2-3k_m+4-mp-3n\le -1.$ Therefore
      \begin{equation}\label{Lemma54Inequality3}
        (k_m-1)^2-mp\le 3n+k_m-4\le 4n.
      \end{equation}

      It is impossible that $(k_m-1)^2-mp=4n,$ otherwise $(k_m-1)^2-1=(m+1)p$, hence $p$ would divide $(k_m-2)k_m$ which forces $k_m=0$ or $k_m=2$, which contradicts $k_m\ge 3.$ Hence $(k_m-1)^2-mp\le 4n-1,$ however, if $(k_m-1)^2-mp=4n-1,$ then $k_m=1$ which is again, impossible. Then $(k_m-1)^2-mp<p$ and hence $\r{(k_m-1)^2}=(k_m-1)^2-mp.$

      Now $\Gamma(k_m)=\r{(k_m-1)^2}+\r{k_m+1-3n}=(k_m-1)^2-mp+k_m+1-3n-p$ \ as \ $k_m\le n+2\le 3n-2.$
      Now Inequality~\ref{Lemma54Inequality1} implies $\Gamma(k_m)=k_m^2-k_m+2-3n-mp+p\ge p,$ hence $k_m$ is a jump.
  \item
      To finish the proof, we observe that from Inequalities~\ref{Lemma54Inequality2} and ~\ref{Lemma54Inequality3} we obtain
      \[3n-k_m\le\r{(k_m-1)^2}=(k_m-1)^2-mp\le3n+k_m-4.\]
\end{enumerate}

\ok

The following lemma allows us to find more jumps based on the ones found in  Lemma~\ref{Lemma54k=Qh-isajump}.

%%%%%%%%%%%%%%%%%%%%%%%%%%%%%%%%%%%%%%%%%%       Lemma             %%%%%%%%%%%%%%%%
\begin{lem}\label{Lemma55jumpsTypeB}
   Consider $k_m$ be the jump in  Lemma~\ref{Lemma54k=Qh-isajump} and $k\le n+2.$ If $k_m<k\le 1+\floor{\sqrt{mp+p-1}}$ then $k$ is a jump and \[ 3n+k-3<\rkm=(k-1)^2-mp.\]
\end{lem}
%%%%%%%%%%%%%%%%%%%%%%%%%%%%%%%%%%%%%%%%%%       Proof of Lemma    %%%%%%%%%%%%%%%%
\noindent {\bf Proof.}
  Clearly \m{k_m<k} implies \m{mp\le (k_m-1)^2<(k-1)^2.} Since $k\le\floor{\sqrt{mp+p-1}},\ \ (k-1)^2\le mp+p-1.$ Hence $\rkm=(k-1)^2-mp.$

  \medskip

  Since $k\le  n+2, \ \r{k+1-3n}=k+1-3n-p.$
  We know that $k$  is strictly above the positive root of $\,x^2-x+2-3n-mp,$ hence  $\,k^2-k+2-3n>mp.$ This implies that
   \begin{align*}
     \gk =& \rkm+\rkl =\kt+p\\
         =&k^2-k+2-3n-mp+p\ge p.
   \end{align*}

  \nn and then $k$ is a jump.

  \medskip

  \nn From Lemma~\ref{Lemma54k=Qh-isajump}, \m{\rpi<k_m\le k-1,} hence $\,\rpi[3]<k.\,$ Therefore $0<k^2-3k+4-mp-3n,$ then \[ 3n+k-3<(k-1)^2-mp=\rkm.\] \ok

%%%%%%%%%%%%%%%%%%%%%%%%%%%%%%%%%%%%%%%%%%       Remark            %%%%%%%%%%%%%%%%
\begin{rem}\label{Observation56}
  Note that by Lemma ~\ref{Lemma54k=Qh-isajump}, the jumps $k_m$ satisty $3n-k_m-1<\rkm\le 3n+k_m-4$ and by Lemma~\ref{Lemma55jumpsTypeB}, the jumps $k>k_m$ satisfy $3n+k-3<\rkm.$ Now if $k$ is not a jump and $2\le k\le n+2$, then necessarily $\rkm<3n-k-1$ as the following lemma shows.
\end{rem}

The following two lemmas are about the total residues and will help us to count the total amount of jumps in the interval $[1,4n].$ We will only prove the fist one as the proof of the second one is similar.

%%%%%%%%%%%%%%%%%%%%%%%%%%%%%%%%%%%%%%%%%%       Lemma             %%%%%%%%%%%%%%%%
\begin{lem}\label{Lemma57ThreeChunks}
 Let $n\ge 2 $ and $2\le k\le n+2$ and $p=4n-1$ prime.
 \begin{enumerate}[(a)\,]
   \item If $\rkm<3n-k-1$ then \, $\gk<p$ \, and\, $\gpk<p.$
   \item If $3n-k-1\le \rkm< 3n+k-3$ then \, $\gpk<p\le \gk.$
   \item If $3n+k-3\le\rkm$ then \, $p\le \gpk$\, and \, $p\le \gk.$
 \end{enumerate}
\end{lem}
%%%%%%%%%%%%%%%%%%%%%%%%%%%%%%%%%%%%%%%%%%       Proof of Lemma    %%%%%%%%%%%%%%%%
%\noindent {\Proof.}
\noindent {\bf Proof.}

 Note that $-(4n-1)\le k+1-3n\le -1,$ therefore $\rkl=k+1-3n+p.$ Similarly, $k\le n+2$ implies $0\le p+3-k-3n\le p-1$ and hence $\rk{p+3-k-3n}=p+3-k-3n.$
 \begin{enumerate}[(a)\,]
   \item
       Clearly
        \begin{align*}
          \gk =& \rkm+\rkl, \\
              =& \rkm+p+k+1-3n<p.\\
          \gpk =& \r{(p-(k-1))^2}+\r{p+3-k-3n}, \\
               =& \rkm+p+3-k-3n<p.
        \end{align*}
   \item Here $\gk=\rkm+p+k+1-3n\ge p$  and
        $\gpk=\rkm+p+3-k-3n<p.$
   \item Finally $\gk=\rkm+p+k+1-3n\ge p$  and
        $\gpk=\rkm+p+3-k-3n\ge p.$
 \end{enumerate}
  Note that $2\le k\le n+2$ implies $0\le p+2-k\le p,$ hence $\gpk$ is defined. Also when $k=0$ or $k=1,\ \ \rk<3n-k-1$ and $\gk<p$. Then $\gpk$ is not defined as $p+2-k>p.$
\ok

%%%%%%%%%%%%%%%%%%%%%%%%%%%%%%%%%%%%%%%%%%       Lemma             %%%%%%%%%%%%%%%%
\begin{lem}\label{Lemma58}
 Let $n\ge 2, p=4n-1$ prime and $n+3\le k\le 2n.$
 \begin{enumerate}[(a)\,]
   \item If $\rkm<k-n-2$ then \, $\gk<p$ \, and\, $\gpk<p.$
   \item If $k-n-2\le \rkm< 3n-k-1$ then $\gk<p\le \gpk.$
   \item If $3n-k-1\le\rkm$ then \, $p\le \gk$\, and \, $\gpk.$
 \end{enumerate}
\end{lem}
%%%%%%%%%%%%%%%%%%%%%%%%%%%%%%%%%%%%%%%%%%       Proof of Lemma    %%%%%%%%%%%%%%%%
\noindent {\bf Proof.}
(very similar to the one of Lemma~\ref{Lemma57ThreeChunks}.)
\ok

\section{Splitting the Jumps}\label{Section:SplittingTheJumps}

\nn For future reference we define the following sets.

%@%%%%%%%%%%%%%%%%%%%%%%%%%%%%%%%%%%%%%%%%%       Notation         %%%%%%%%%%%%%%%%

\begin{nota}\label{NotationA'sAndB's} For $p=4n-1$ prime, denote
  \[
   \begin{split}
     \hAle   & = \left\{k\in \{2,...,n+2\}\ :\ \rkm<3n-k-1 \, \right\}, \\
     \hAbet  & = \left\{k\in \{2,...,n+2\}\ :\ \rkm\in [3n-k-1,3n+k-3) \, \right\}, \\
     \hAge   & = \left\{k\in \{2,...,n+2\}\ : \ \rkm\ge 3n+k-3  \, \right\}, \\[2mm]
     \hBle   & = \left\{k\in \{n+3,...,2n\}\ : \ \rkm<k-n-2\, \right\}, \\
     \hBbet  & = \left\{k\in \{n+3,...,2n\}\ : \ \rkm\in [ k-n-2,3n-k-1) \,\right\},\\
     \hBge   & = \left\{k\in \{n+3,...,2n\}\ : \ \rkm\ge 3n-k-1 \,\right\}. \\
   \end{split}
  \]
\end{nota}

%%%%%%%%%%%%%%%%%%%%%%%%%%%%%%%%%%%%%%%%%%       Lemma             %%%%%%%%%%%%%%%%
\begin{lem}\label{Lemma60TwoinA<TwoinA[)}
  Let $\ell\in\Z$ and $p=4n-1$ prime.
  \begin{enumerate}[(a)\,]
    \item Let $\ell\ge 6$. If $n=4\ell+2$ or $n=4\ell+3$  then $n-3,n-1,n+1\in\hAbet$ and $n-2,n,n+2\in\hAle.$
    \item Let $\ell\ge 4$. If $n=4\ell+1$ or $n=4\ell+4$  then $n-3,n-1,n+1\in\hAle$ and $n-2,n,n+2\in\hAbet.$
  \end{enumerate}
\end{lem}
%%%%%%%%%%%%%%%%%%%%%%%%%%%%%%%%%%%%%%%%%%       Proof of Lemma    %%%%%%%%%%%%%%%%
\noindent {\bf Proof.}
  Table~\ref{TableResidues} summarizes the residues of $(k-1)^2$ for $k=n-3,n-2,n-1,n,n+1$ and $n+2$ given the four different cases for $n$ which are $4\ell+1,4\ell+2,4\ell+3$ and $4\ell+4,\ \ \ell\ge 4.$

  We only verify the first column of Table~\ref{TableResidues}, that is, we will compute the residues of $(k-1)^2$ when $k=n-3.$
  \begin{eqnarray*}
  % \nonumber to remove numbering (before each equation)
    (4\ell-3)^2 &=& (16\ell+3)(\ell-2)+5\ell+15, \ \ 0\le 5\ell+15 \le 16\ell+2,  \\
    (4\ell-2)^2 &=& (16\ell+7)(\ell-2)+9\ell+18, \ \ 0\le 9\ell+18 \le 16\ell+6, \\
    (4\ell-1)^2 &=& (16\ell+11)(\ell-2)+13\ell+23, \ \ 0\le 13\ell+23 \le 16\ell+10, \\
    (4\ell)^2 &=& (16\ell+15)(\ell-1)+\ell+15, \ \ 0\le \ell+15 \le 16\ell+14.
  \end{eqnarray*}
 \begin{table}[h]\small
   \[
   \begin{tblr}{
    hlines = {black},
    vlines = {black},
    row{2} = {c},
    cell{3,5}{1,2,3,4,5,6,7} = {gray9,m},
    cell{2,4,6}{1,2,3,4,5,6,7} = {m},
    cell{1}{2} = {r=1,c=6}{c},
   }
     & k &&&&&\\
           & n-3             & n-2             & n-1        & n & n+1 & n+2 \\
    {n=4\ell+1\\ \ns p=16\ell+3} & 5\ell+15  & 13\ell+10  &5\ell+4 & 13\ell+3 & 5\ell+1 & 13\ell+4 \\
    {n=4\ell+2\\ \ns p=16\ell+7} & 9\ell+18  & \ell+8  &9\ell+7 & \ell+1 & 9\ell+4 & \ell+2 \\
    {n=4\ell+3\\ \ns p=16\ell+11} & 13\ell+23 & 5\ell+11  &13\ell+12 & 5\ell+4 & 13\ell+9 & 5\ell+5 \\
    {n=4\ell+4\\ \ns p=16\ell+15} & \ell+15  & 9\ell+16  &\ell+4 & 9\ell+9 & \ell+1& 9\ell+10 \\
    \end{tblr}
    \]
      \caption{Residues of $(k-1)^2$ when $k=n-3,n-2,...,n+2$ for the different cases of $n.$}
      \label{TableResidues}
  \end{table}
  For a given $n,k$, denote $I_n^k=[3n-k-1,3n+k-3)$ and $\Delta_n^k=\rkm.$
  \begin{enumerate}[(a)\,]
    \item Consider $n=4\ell+2$ and $k=n-3$ then $I_n^k=[8\ell+6,16\ell+2).$ According to Table~\ref{TableResidues}, $\Delta_n^k=9\ell+18.$ We observe that if $\ell\ge 3,\ \Delta_n^k\in I_n^k.$  Similarly, if $k=n-1, I_n^k=[8\ell+4,16\ell+4)$ and $\Delta_n^k=9\ell+7\in I_n^k$ for $\ell\ge 1.$ If $k=n+1, I_n^k=[8\ell+2,16\ell+6)$ and $\Delta_n^k=9\ell+4\in I_n^k$ for $\ell\ge 0.$ Hence for $\ell\ge 6,\ \ n-3, n-1, n+1\in \hAbet.$

        Notice that if $k=n-2, \Delta_n^k=\ell+8\notin [8\ell+5,16\ell+3)= I_n^k$ when $\ell\ge 1.$ Likewise, if $k=n, \Delta_n^k=\ell+1\notin [8\ell+3,16\ell+5)= I_n^k$ when $\ell\ge 0.$ Finally, if $k=n+2, \Delta_n^k=\ell+2\notin [8\ell+1,16\ell+7)= I_n^k$ when $\ell\ge 1.$ Therefore $n-2,n,n+2\notin \hAbet$ when $\ell\ge 1,$ in fact, $n-2,n,n+2\in \hAle$ for $\ell\ge 1.$ The case $n=4\ell+3$ is analogous.
    \item This case is done analogously.
  \end{enumerate}
\ok

%%%%%%%%%%%%%%%%%%%%%%%%%%%%%%%%%%%%%%%%%%       Observation 60.5          %%%%%%%%%%%%%%%%
\begin{rem}\label{Observation60.5}
If $n>3$ either $\gm[n+2]<p\le\gm$ or $\gm<p\le\gm[n+2].$
\end{rem}
%%%%%%%%%%%%%%%%%%%%%%%%%%%%%%%%%%%%%%%%%%       Proof of Observation    %%%%%%%%%%%%%%%%
\noindent {\bf Proof.}
From Lemma~\ref{Lemma60TwoinA<TwoinA[)}, if $n=4\ell+2$ or $n=4\ell+3$ and $\ell\ge 6$, then $n+1\in\hAbet$ and $n+2\in\hAle.$ By Lemma~\ref{Lemma57ThreeChunks},  $\gm[n+2]<p\le\gm$. From Lemma~\ref{Lemma60TwoinA<TwoinA[)}, if $n=4\ell+1$ or $n=4\ell+4$ and $\ell\ge 4$, then $n+1\in\hAle$ and $n+2\in\hAbet.$ By Lemma~\ref{Lemma57ThreeChunks},  $\gm<p\le\gm[n+2]$.

We only need to verify the cases $n=5,6,8,11,12,15$ and $18$ as the cases $n=4,7,9,10,13,16$ and $19$ do not give prime numbers.
\begin{enumerate}[(i)]
  \item If $n=5$ then $n+2 \in\hAbet=\{5,7\}$ and $n+1 \in\hAle=\{2, 3, 4, 6\}.$
  \item If $n=6$ then $n+1 \in\hAbet=\{5,7\}$ and $n+2 \in\hAle=\{2, 3, 4, 6, 8\}.$
  \item If $n=8$ then $n+2 \in\hAbet=\{6, 8, 10\}$ and $n+1 \in\hAle=\{2, 3, 4, 5, 7, 9\}.$
  \item If $n=11$ then $n+1 \in\hAbet=\{7,10,12\}$ and $n+2 \in\hAle.$
  \item If $n=12$ then $n+2 \in\hAbet=\{7, 10, 12, 14\}$ and $n+1 \in\hAle.$
  \item If $n=15$ then $n+1 \in\hAbet=\{8,11,14,16\}$ and $n+2 \in\hAle.$
  \item If $n=18$ then $n+1 \in\hAbet=\{8, 12, 15, 17, 19\}$ and $n+2 \in\hAle.$
\end{enumerate}
An application of Lemma~\ref{Lemma57ThreeChunks} gives us the result.

Observe that when $n=3,\hAbet=\{4,5\},\hAle=\{2,3\}$, in this case both $\gm=16$\, and \,$\gm[n+2]=13$ are greater than $p=11.$  \ok

The following two lemmas allow us to compute specifically the cardinality of $\hAge.$
%%%%%%%%%%%%%%%%%%%%%%%%%%%%%%%%%%%%%%%%%%       Lemma             %%%%%%%%%%%%%%%%
\begin{lem}\label{Lemma58.5}
 Let $n> 3,\ p=4n-1$ prime and $M_0=\floor{(n^2-4n+5)/p}$. Define $\ell_m=1+\floor{\sqrt{mp+p-1}}$ and consider $k_m$ defined as in Lemma~\ref{Lemma54k=Qh-isajump}. Let $2\le k\le n+2$ and $0\le m \le M_0-1.$

  If $k<k_0, \, k>\ell_{M_0}$ or $\ell_m<k\le k_{m+1}$ then $k\notin \hAge.$ Also if $k<k_0$ or $\ell_m<k<k_{m+1}$ then $k$ is not a jump.
\end{lem}
%%%%%%%%%%%%%%%%%%%%%%%%%%%%%%%%%%%%%%%%%%       Proof of Lemma    %%%%%%%%%%%%%%%%
\noindent {\bf Proof.}
Consider $k<k_0=1+\floor{(1+\sqrt{12n-7})/2}$. Then $(2k-1)^2<12n-7$ from which $(k-1)^2<3n-k-1\le 4n-1.$ Hence $\rkm=(k-1)^2<3n-k-1.$ Therefore by Lemma~\ref{Lemma57ThreeChunks}, $k\in\hAle$ and $k$ is not a jump.

\medskip

\nn Consider $k>\ell_{M_0}.$ Notice that $\ell_{M_0}=1+\floor{\sqrt{Mp}}\ge n-1.$ By Lemma~\ref{Lemma60TwoinA<TwoinA[)}, either $k\in\hAbet$ or $k\in\hAle.$ If $k=k_{m+1},$ by Lemma~\ref{Lemma54k=Qh-isajump} (iv), $k\in \hAbet.$

\medskip

\nn Finally consider $\ell_m<k<k_{m+1}.$  Then \[1+\sqrt{(m+1)p}<k<\yp[1]{(m+1)}.\]

Hence
\[(m+1)p<(k-1)^2<(m+1)p+3n-k-1.\]

Therefore $\rkm=(k-1)^2-(m+1)p<3n-k-1.$ By Lemmas~\ref{Lemma57ThreeChunks} and~\ref{Lemma60TwoinA<TwoinA[)}, $k$ is not a jump and $k\in\hAle.$ \ok

%%%%%%%%%%%%%%%%%%%%%%%%%%%%%%%%%%%%%%%%%%       Lemma 53             %%%%%%%%%%%%%%%%
\begin{lem}\label{Lemma72}
  Consider $n, p$ and $M_0$ as in Lemma~\ref{Lemma54k=Qh-isajump}. Then $k\in\hAge$ if and only if there is $m\in\{0,1,...,M_0\}$ such that $k_m<k\le 1+\floor{\sqrt{mp+p-1}}.$
\end{lem}
%%%%%%%%%%%%%%%%%%%%%%%%%%%%%%%%%%%%%%%%%%       Proof of Lemma    %%%%%%%%%%%%%%%%
\noindent {\bf Proof.}

From Lemma~\ref{Lemma55jumpsTypeB}, if $k_m<k\le\ell_m$ then $k$ is a jump and $\rkm>3n+k-3.$ Note that $2\le k_0<k.$  Since $m\le M_0$ we have $\ell_m\le n+2,$ hence $2\le k \le n+2.$ Therefore $k\in\hAge.$

Conversely, let $k\in\hAge.$ Then $2\le k \le n+2$ and $\rkm\ge 3n+k-3.$ Observe that
{\small
\[ [2,n+2]=[2,k_0]\cup(k_0,\ell_0]\cup(\ell_0,k_{1}]\cup
(k_1,\ell_1]\cup\cdots\cup(k_{M_0},\ell_{M_0}]\cup(\ell_{M_0},n+2].\]
}

If $k<k_0, \, k>\ell_{M_0}$ or $\ell_m<k\le k_{m+1}$, by Lemma~\ref{Lemma58.5}, $k\notin \hAge.$

\medskip

This forces $k$ to be in  $(k_m,\ell_m]$ for some $m\in\{0,...,M_0\},$ which is what we wanted. \ok

\section{Counting the Jumps}\label{Section:CountingTheJumps}

In this section we will relate the jumps in the different sets $\hAle,\hAbet,\hAge,\hBle,\hBbet$ and $\hBge.$  We will compute different cardinalities  when possible. In this section we consider $n>3$ and $p=4n-1$ a prime number.

%%%%%%%%%%%%%%%%%%%%%%%%%%%%%%%%%%%%%%%%%%       Theorem           %%%%%%%%%%%%%%%%
\begin{thm}\label{Lemma59TheoremA=>Bijective}
  The function $f:\hAge\too \hBle$ defined by
  \[f(k)=2n+2-k,\]
  is well-defined and bijective.
\end{thm}
%%%%%%%%%%%%%%%%%%%%%%%%%%%%%%%%%%%%%%%%%%       Proof of Theorem    %%%%%%%%%%%%%%
\noindent {\bf Proof.}

Consider $k\in\hAge.$ By Lemmas~\ref{Lemma55jumpsTypeB},~\ref{Lemma72} and Observation~\ref{Observation56}, there is $m\in\{0,1,...,M_0\}$ with $M_0=\floor{(n^2-4n+5)/p}$ such that $k_m<k\le 1+\floor{\sqrt{mp+p-1}}$ and $3n+k-3\le (k-1)^2-mp=\rkm,$  where $k_m$ is the jump given in Lemma~\ref{Lemma57ThreeChunks}. Hence
  \begin{eqnarray}
     0\le k^2-3k-3n+4-mp, \label{59-1}\\
     k^2-2k+1-mp\le 4n-2. \label{59-2}
  \end{eqnarray}

Let $k_f=f(k)=2n+2-k.$ We will first prove that $k_f\in \hBle.$ Now
   \[ (k_f-1)^2=(n-k+m+2)p+\keh. \]
Let $w=\keh.$ From Inequality~\ref{59-1}, $w\ge -1.$ From Remark~\ref{Note53.5k^2-3k+4-3n=hp}, $k^2-3k+4-3n-mp\neq 0,$ hence $w\ge 0.$

Note that in the case of equality in Inequality~\ref{59-2}, we would have $(k-1)^2=mp+p-1,$ hence the congruency $x^2\equiv -1 \pmod p$ would have a solution, which is impossible. Therefore
      \begin{equation}\label{59-3}
        k^2-2k+1-mp<4n-2.
      \end{equation}
From Inequality~\ref{59-2},
      \begin{equation}\label{59-4}
        w=k^2-2k+1-mp+2-k-3n<n-k.
      \end{equation}
Clearly $n-k\le 4n-2,$ hence $w<p-1.$ Therefore $\rkm[k_f]=w.$
From Inequality~\ref{59-4}, $0\le w<n-k,$ this forces $k\le n-1$. Also from Inequality~\ref{59-4}, since $k_f-n-2=n-k,$ we conclude that $\rkm[k_f]=w<k_f-n-2.$

Finally, $2\le k\le n-1$ implies $n+3\le k_f\le 2n,$ hence $f$ is well defined as $k_f\in\hBle.$

Clearly $f$ is injective. Take now $\wk\in\hBle.$ \ Then $n+3\le\wk\le 2n$ and
      \begin{equation}\label{59-5}
         \rkm[\wk]<\wk-n-2.
      \end{equation}

      Consider $k=2n+2-\wk$ and $\wh\in\Z$ with
      \begin{equation}\label{59-6}
         \wh\le(\wk-1)^2<\wh\cdot p+p.
      \end{equation}
Then $2\le k\le n-1$ and
\[(k-1)^2=(n-\wk+\wh)\cdot p+\khhh.\]
Consider $m=n-\wk+\wh$ and $u=\khhh.$

Since $\wk\le 5n,$ Inequality~\ref{59-6} implies
      \begin{equation}\label{59-7}
           0\le n-\wk+1-p\le \khkh.
      \end{equation}

Also from Inequality~\ref{59-5}, we have  $\khkh<p-1,$ therefore $0\le u<p-1.$ Hence $\rkm=(k-1)^2-mp=u.$

Finally Inequality~\ref{59-7} implies $u\ge 5n-\wk\ge5n-1-\wk=3n+k-3.$
Therefore $\rkm\ge 3n+k-3,$ i.e. $k\in\hAge.$ Clearly $f(k)=\wk$, then $f$ is bijective.\ok

%%%%%%%%%%%%%%%%%%%%%%%%%%%%%%%%%%%%%%%%%%       Theorem           %%%%%%%%%%%%%%%%
\begin{lem}\label{Lemma61hSatisfies}
  Let $y,z$ be the last two elements in $\hAbet.$ If $k\in \hAbet-\{y,z\}$ then
      \begin{equation}\label{61-1}
           mp+p\le (n-2)^2+1,
      \end{equation}
  where $m=m(k)=\floor{(k-1)^2/p}.$
\end{lem}
%%%%%%%%%%%%%%%%%%%%%%%%%%%%%%%%%%%%%%%%%%       Proof of Theorem    %%%%%%%%%%%%%%
\noindent {\bf Proof.}
 \begin{enumerate}[(a)]
   \item Case $n=4\ell+2$ or $n=4\ell+3, \ \ell\ge 6.$ By Lemma~\ref{Lemma60TwoinA<TwoinA[)}, $k\in\hAbet-\{y,z\}$ implies $k\le n-3.$ Hence $m\le \ell-2=\floor{(n-4)^2/p}.$
       If $n=4\ell+2$ then $mp+p\le(\ell-1)(16\ell+7)\le16\ell^2+1=(n-2)^2+1.$

       If $n=4\ell+3$ then $mp+p\le(\ell-1)(16\ell+11)\le16\ell^2+8\ell+2=(n-2)^2+1$.
   \item Case $n=4\ell+1$ or $n=4\ell+4, \ \ell\ge 4.$ By Lemma~\ref{Lemma60TwoinA<TwoinA[)}, $k\in\hAbet-\{y,z\}$ implies $k\le n-2.$
       If $n=4\ell+1$ then $m\le \floor{(n-3)^2/p}=\ell-2.$ Hence $mp+p\le(\ell-1)(16\ell+3)\le 16\ell^2-8\ell+2=(n-2)^2+1.$

       If $n=4\ell+4$ then $m\le \floor{(n-3)^2/p}=\ell-1.$ Hence $mp+p\le\ell\cdot(16\ell+15)\le 16\ell^2+16\ell+5=(n-2)^2+1.$
   \item The only left cases are $n=5,6,8,11,12,15$ and 18 as the choices $n=4, 7, 9, 10 ,1 3, 14 ,1 6, 19, 22, 23$ do  not provide prime numbers.

       \begin{enumerate}[(i)]
         \item If $n=5$ or $6$ then $\hAbet$ has only two elements. Hence Inequality~\ref{61-1} is trivial.
         \item If $n=8$ then $\hAbet=\{6, 8, 10\}$ and $k=6.$ Therefore $m(k)=0=\floor{(k-1)^2/p}=\floor{25/31}$ clearly satisfies Inequality~\ref{61-1}.
         \item If $n=11$ then $\hAbet=\{7, 10, 12\}$ and $k=7.$ Therefore $m(k)=0$ clearly satisfies Inequality~\ref{61-1}.
         \item If $n=12$ then $\hAbet=\{7, 10, 12, 14\}.$ If $k=10,$ then $m(k)=1=\floor{(k-1)^2/p}=\floor{81/47}$ clearly satisfies Inequality~\ref{61-1} as $94\le 101$. Clearly $m(7)$ also does.
         \item If $n=15$ then $\hAbet=\{8,11,14,16\}.$  If $k=11,$ then $m=1=\floor{(k-1)^2/p}=\floor{100/59}$ clearly satisfies Inequality~\ref{61-1} as $118\le 169$. Clearly $m(8)$ also does.
         \item If $n=18$ then $n+1 \in\hAbet=\{8, 12, 15, 17, 19\}.$  If $k=15,$ then $m(k)=2=\floor{(k-1)^2/p}=\floor{196/71}$ clearly satisfies Inequality~\ref{61-1} as $213\le 256$. Clearly $m(8), m(12)$ also  satisfy Inequality~\ref{61-1}.
       \end{enumerate}\ok
 \end{enumerate}

%%%%%%%%%%%%%%%%%%%%%%%%%%%%%%%%%%%%%%%%%%       Lemma             %%%%%%%%%%%%%%%%
\begin{lem}\label{Lemma62Bbet}
  $k\in\hBbet$ if and only if there is an integer $m$ with $1\le m\le \floor{(n^2-4n+5)/p}$ and
        \[k=2n-\floor{\sqrt{mp-1}}.\]
\end{lem}
%%%%%%%%%%%%%%%%%%%%%%%%%%%%%%%%%%%%%%%%%%       Proof of Lemma    %%%%%%%%%%%%%%%%
\noindent {\bf Proof.}
\begin{description}
  \item[$\Rightarrow$)]  Let $k\in\hBbet$ and define $\alpha=2n+1-k.$ Hence $1\le \alpha\le n-2$ and \[(k-1)^2=4n^2+\alpha^2-4n\alpha=(n-\alpha+m)p+\alpha^2+n-\alpha-mp,\]
  where $m$ satisfies $mp\le \alpha^2+n-\alpha<(m+1)p.$ Therefore
  $\rkm=\alpha^2+n-\alpha-mp.$ \ Since $k-n-2\le \rkm<3n-k-1,\ \ mp-1\le\alpha^2$ and $(\alpha-1)^2<mp-1.$

  Hence $\sqrt{mp-1}\le\alpha<\sqrt{mp-1}+1.$ Since $x^2\equiv -1 \pmod p$ has no solution, $\sqrt{mp-1}$ is not an integer. Therefore $\alpha=\floor{\sqrt{mp-1}}+1.$

  Then $k=2n-\floor{\sqrt{mp-1}}$.

  Since $0\le\alpha^2+n-\alpha, \ m\ge 0.$ Now $(\alpha-1)^2<mp-1$ implies $1\le m.$ Clearly $mp\le\alpha^2+1\le(n-2)^2+1.$ Therefore $m\le (n^2-4n+5)/p.$

  \item[$\Leftarrow$)] Consider an integer $m$ with  $1\le m\le (n^2-4n+5)/p$ and
  $k=2n-\floor{\sqrt{mp-1}}.$

  Since $mp-1\le (n-2)^2, \ \  n+2\le k.$
  Clearly $n+2=k$ leads us to an integer solution of $x^2\equiv -1 \pmod p,$ therefore  $n+3\le k.$ Also $1\le m$ implies  $k\le 2n.$

  Consider $\alpha=\floor{\sqrt{mp-1}}+1.$ Then $\sqrt{mp-1}\le \alpha<\sqrt{mp-1}+1.$ Clearly $\alpha\neq \sqrt{mp-1}+1$ (otherwise $x^2\equiv -1 \pmod p$ has an integer solution), then $\sqrt{mp-1}< \alpha < \sqrt{mp-1}+1$.

  Hence $mp-1<\alpha^2$ and $(\alpha-1)^2<mp-1.$ Since $\alpha<n-1, \ mp\le\alpha^2\le \alpha^2+n-\alpha<n+\alpha-2+mp<mp+p.$ Since
  \[(k-1)^2=(2n-\alpha)^2=(n-\alpha+m)p+\alpha^2+n-\alpha-mp,\]
  $\rkm=\alpha^2+n-\alpha-mp.$  Finally, from $k-n-2=n-\alpha-1<\rkm<n+\alpha-2=3n-k-1,$ we conclude that $k\in\hBbet.$ \ok
\end{description}

%%%%%%%%%%%%%%%%%%%%%%%%%%%%%%%%%%%%%%%%%%       Theorem           %%%%%%%%%%%%%%%%
\begin{thm}\label{Theorem63TheoremAMinusTwo=>Bijective}
  Let $y,z$ the last two elements of $\hAbet.$  For  $k\in \hAbet-\{y,z\}$, consider $m=\floor{(k-1)^2/p}$ and $u_0=u_0(k)$ the first integer less than or equal to $n+2$ such that $x=u_0$ satisfies

      \begin{equation}\label{63-1}
         (m+1)\,p\le(k+x-1)^2+1.
      \end{equation}

  Then, the function $f:\hAbet-\{y,z\}\too \hBbet$ defined by
  \[f(k)=k_f=2n+2-u_0-k,\]
  is well-defined and bijective.
\end{thm}
%%%%%%%%%%%%%%%%%%%%%%%%%%%%%%%%%%%%%%%%%%       Proof of Theorem    %%%%%%%%%%%%%%
\noindent {\bf Proof.}
  First, we will prove that $f$ is well-defined. It is not hard to check that $u_0=\floor{\sqrt{(m+1)p}}+2-k$. The definition of $u_0$ implies that $u_0,u_0+1,...$ satisfy Inequality~\ref{63-1} but $u_0-1$ does not. Hence $u_0$ satisfies
  \begin{eqnarray}
         4n+mp\le k^2+2k(u_0-1)+u_0^2-2u_0+3,\label{63-2}\\
         k^2+2k(u_0-2)+u_0^2-4u_0+6< 4n+mp. \label{63-3}
  \end{eqnarray}
  If $k^2=mp+3n+3k-2$ for $2\le k\le n+2$ and we define $x_0=2n-k+1,$ then $x_0^2=(n+m-k+2)p+1.$ Therefore $p$ divides $(x_0-1)(x_0+1)$, however since $n\ge 3$ and $2\le k \le n+2,$ we have that
  \[1\le 2n-k =x_0-1<x_0+1=2n+2-k\le 4n-2,\]
  which is impossible. Therefore
      \begin{equation}\label{63-4}
         k^2\neq mp+3n+3k-2.
      \end{equation}
  Observe that $u_0\ge 1$ as $x=0$ does not satisfy Inequality~\ref{63-1}. Also $u_0\le n+2$ as $(k+n+1)^2\ge(k-1)^2+(n+2)^2\ge mp+p-1.$

  \nn To shorten notation, define
   \[ \left\{\at[lllllllllllll]{m_f & =& n+2-k-u_0+m,\\ \Delta & = &\rkm,\\ \Delta_f & = &\rkm[k_f],\\ w_f & = & \Delta+k(2u_0-1)-p+n+u_0^2-3u_0+1. }\right. \]

  To check that $f$ is well-defined, we need to verify that $n+3\le k_f\le 2n$ and $k_f-n-2\le\Delta_f<3n-k_f-1.$ It is not hard to check that
     \[(k_f-1)^2=m_f\cdot p+w_f.\]
  Notice that $\Delta=\rkm=(k-1)^2-mp\le 3n-k-1.$ From Inequality~\ref{63-4},
  $\Delta\ge 3n-k.$ Therefore,
 \begin{eqnarray*}
      w_f&\ge & 3n-k+ +k(2u_0-1)-p+n+u_0^2-3u_0+1, \\
      &=&k(2u_0-2)+u_0^2-3u_0+2,\\
      &\ge& u_0^2-3u_0+2=(u_0-2)(u_0-1)\ge 0.
 \end{eqnarray*}
 From Inequality~\ref{63-3},
 \begin{eqnarray*}
      w_f&= & k^2+2k(u_0-2)+u_0^2-4u_0+6-4n-mp+k+u_0+n-3, \\
      &<&k+u_0+n-3\le 4n-1.
 \end{eqnarray*}

 This shows that $w_f=\rkm[k_f].$ Since $3n-k_f-1=n-3+u_0+k,\ w_f<3n-k_f-1.$  Since $k_f-n-2=n-u_0-k,$ Inequality~\ref{63-2} implies that $w_f\ge k_f-n-2.$ Therefore $k_f-n-2\le\rkm[k_f]<3n-k_f-1.$ By Lemma~\ref{Lemma61hSatisfies}, $x=n-k-1$ satisfies Inequality~\ref{63-1}, hence $1\le u_0\le n-k-1.$ Therefore $n+3\le k_f\le 2n.$ This proves that $k_f\in\hBbet$ and thus $f$ is well defined.

 \bigskip

 \nn Consider $ k_0$ and $k_1$ such that $k_f=f(k_0)=f(k_1).$ \ Take $u_0,u_1,m_0$ and $m_1$ such that $m_0=\floor{(k_0-1)^2/p}, m_1=\floor{(k_1-1)^2/p}$ and $u_0, u_1$ are the first integers such that Inequality~\ref{63-1} holds with $k=k_0$ and $k=k_1$ respectively.

 Now $f(k_0)=f(k_1)$ implies that $u_0+k_0=u_1+k_1$. Since
    \[m_f=n+2-k_0-u_0+m_0=n+2-k_1-u_1+m_1,\]
 we conclude $m_0=m_1.$ Since
    \[w_f=(k_0+u_0-1)^2-u_0-m_0p-p+n+1=(k_1+u_1-1)^2-u_1-m_1p-p+n+1,\]
 $u_0=u_1$ and consequently $k_0=k_1.$ Therefore $f$ is injective.

 \bigskip

  \nn Take now $k_f\in\hBbet.$ Let $m_f=\floor{(k_f-1)^2/p}$. Then $k_f-n-2\le \Delta_f<3n-k_f-1.$ Notice that $x=0$ satisfies the inequality
      \begin{equation}\label{63-6}
         0\le (k_f+x-1)^2+1- m_fp-px.
      \end{equation}
Let $v_0$ the first integer greater than or equal to 1 such that $x=v_0$ does not satisfy Inequality~\ref{63-6}. Hence $v_0=x$ satisfies the equivalent inequalities
  \begin{eqnarray}
        (k_f+x-1)^2+1- m_fp-px<0,\nonumber\\[-.25cm]
        \label{63-7}\\[-.25cm]
         \Delta_f+2k_fx+(x-1)^2-px< 0. \nonumber
  \end{eqnarray}

Consider $k=2n+2-v_0-k_f,\ m=n-k_f+m_f$ and $w=\Delta_f+k_f(2v_0-1)+n+(v_0^2-3v_0+1)-v_0p+p.$ Then
  \begin{eqnarray*}
         & \hspace{-.5cm} (k-1)^2 &\\
         &=& \hspace{-.25cm}(n+1-k_f+m_f)p+k_f^2+k_f(2v_0-3)+n+2+v_0^2-3v_0-v_0p-mp, \\
         &=& \hspace{-.25cm}(n-k_f+m_f)p+\Delta_f+k_f(2v_0-1)+n+(v_0^2-3v_0+1)-v_0p+p, \\
         &=& \hspace{-.25cm}mp+w.
  \end{eqnarray*}
By Lemma~\ref{Lemma62Bbet}, there is an integer $m, 1\le m\le (n^2-4n+5)/p$ such that $k_f=2n-\floor{\sqrt{mp-1}}.$ \ Also from the proof of Lemma~\ref{Lemma62Bbet} if $\alpha=2n-k_f+1$ then $m_f=n-\alpha+m.$  Take $x_0=\alpha-1=2n-k_f.$ Notice that $x_0\ge 1$ as $\sqrt{mp-1}\ge 1.$\ Substituting $x=x_0$ into $(k_f+x-1)^2+1- m_fp-px$ gives us
\[(2n-1)^2+1-(n-\alpha+m)p-p(\alpha-1)=n+3-mp.\]
Since $m\ge 1,\ n+3-mp<0.$ Then $x=x_0$ satisfies Inequality~\ref{63-7}, therefore $v_0$ exists and $v_0\le 2n-k_f.$  Since $x=v_0-1$ satisfies Inequality~\ref{63-6}, we have that
\[0\le\Delta_f+k_f(2v_0-2)+(v_0-2)^2+p-v_0p.\]
Hence $w\ge k_f+n+v_0-3=3n-k-1.$ Since $k_f\ge n+3$ and $v_0\ge 1$ we conclude $w\ge 0$ and $2\le k\le n-2.$

\nn Since $v_0$ satisfies Inequality~\ref{63-7},
\[w<n-k_f-v_0+p=5n-k_f-v_0-1=3n+k-3\le 4n-1.\]
Thus $0\le w<p$ and $3n-k-1\le w<3n+k-3.$ This implies that $w=\rkm,\  m=\floor{(k-1)^2/p}$ and hence $k\in \hAbet.$
Lemma~\ref{Lemma60TwoinA<TwoinA[)} and the proof of Lemma~\ref{Lemma61hSatisfies} implies that $\{y,z\}$ is a subset of $\{n-1,n,n+1,n+2\},$ then $k\in\hAbet-\{y,z\}.$

\medskip

\nn Now we will find $u_0=u_0(k).$  Since $\Delta_f\ge k_f-n-2,$
 \begin{eqnarray*}
      (k+v_0-1)^2+1 & = & (2n+1-k_f)^2+1,\\
      & = & (n-k_f+m_f+1)p+\Delta_f+n-k_f+2,\\
      & = & (m+1)p+\Delta_f+n-k_f+2\ge (m+1)p.\\
 \end{eqnarray*}

\nn From $\Delta_f< 3n-k_f-1=p-n-k_f,$ we obtain
 \begin{eqnarray*}
      (k+v_0-2)^2+1 & = & (2n-k_f)^2+1, \\
      & = & (n-k_f+m_f+1)p+\Delta_f+k_f+n-p,\\
      & = & (m+1)p+\Delta_f+k_f+n-p<(m+1)p.\\
 \end{eqnarray*}

 Therefore $u_0=v_0$ and $f(k)=k_f.$ Then $f$ is surjective and hence bijective. \ok

%%%%%%%%%%%%%%%%%%%%%%%%%%%%%%%%%%%%%%%%%%       Corolario         %%%%%%%%%%%%%%%%
\begin{cor}\label{Corollary64|B<|=|A>=|}
  $\big|\hBbet\big|=\floor{(n^2-4n+5)/2}$ and
  \[\big|\hAge\big|=\big|\hBle\big|, \ \ \big|\hAbet\big|=\big|\hBbet\big|+2,\ \ \big|\hAle\big|=\big|\hBge\big|+1.\]
\end{cor}
%%%%%%%%%%%%%%%%%%%%%%%%%%%%%%%%%%%%%%%%%%       Proof of Corollary  %%%%%%%%%%%%%%
\noindent {\bf Proof.}
From Lemma~\ref{Lemma59TheoremA=>Bijective}, $\big|\hAge\big|=\big|\hBle\big|.$ From Theorem~\ref{Theorem63TheoremAMinusTwo=>Bijective}, $\big|\hAbet\big|=\big|\hBbet\big|+2.$ Since $n+1=\big|\hAle\big|+\big|\hAbet\big|+\big|\hAge\big|$ and $n-2=\big|\hBle\big|+\big|\hBbet\big|+\big|\hBge\big|$ we have
 \[n+1=\big|\hAle\big|+\big|\hBle\big|+\big|\hBbet\big|+2=\big|\hAle\big|+n-\big|\hBge\big|,\]
 hence $\big|\hAle\big|=\big|\hBge\big|+1.$
To see that $\big|\hBbet\big|=\floor{(n^2-4n+5)/2}$ it is enough to see that all the $k's$ in Lemma~\ref{Lemma62Bbet} given by each $m$ are all different, which is the case as
\[\sqrt{mp-1}+1<\sqrt{(m+1)p-1}.\]
\ok

%%%%%%%%%%%%%%%%%%%%%%%%%%%%%%%%%%%%%%%%%%       Corolario         %%%%%%%%%%%%%%%%
\begin{thm}\label{Theorem65}
  Under the hypothesis of Theorem~\ref{Theorem63TheoremAMinusTwo=>Bijective},
  \[\Big|\{k\in\Z\ | \ 2\le k\le 4n-1, \ \gk\ge p\}\Big|=2n-2.\]
\end{thm}
%%%%%%%%%%%%%%%%%%%%%%%%%%%%%%%%%%%%%%%%%%       Proof of Corollary  %%%%%%%%%%%%%%
\noindent {\bf Proof.}
Let $J_{\Gamma}=\left\{k\in\Z\ :\ 2\le k\le 4n-2,\ \ \gk\ge p\right\}$. By Lemmas~\ref{Lemma57ThreeChunks} and ~\ref{Lemma58},
\begin{eqnarray*}
    J_{\Gamma}\cap [2,n+2]                  &=&  \hAbet\cup \hAge,\\
    J_{\Gamma}\cap [n+3,2n]                 &=&  \hBge.
\end{eqnarray*}

\begin{eqnarray*}
    \big|J_{\Gamma}\cap [2n+1,3n-2]\big|  && \\
      &&\hspace{-3cm}= \Big|\{2n+1\le k\le 3n-2\ \ | \ \ \gk\ge p\}\Big|, \\
      &&\hspace{-3cm}= \Big|\{2n+1\le p+2-k\le 3n-2\ \ | \ \ \gm[p+2-k]\ge p\}\Big|,  \\
      &&\hspace{-3cm}= \Big|\{n+3\le k\le 2n\ \ |\  \ \gm[p+2-k]\ge p\}\Big| =  \big|\hBbet\big|+\big|\hBge\big|.
\end{eqnarray*}

\begin{eqnarray*}
    \big|J_{\Gamma}\cap [3n-1,4n-1]\big|  && \\
      &&\hspace{-3cm}= \Big|\{3n-1\le k\le 4n-1\ \ | \ \ \gk\ge p\}\Big|, \\
      &&\hspace{-3cm}= \Big|\{3n-1\le p+2-k\le 4n-1\ \ | \ \ \gm[p+2-k]\ge p\}\Big|,  \\
      &&\hspace{-3cm}= \Big|\{2\le k\le n+2\ \ |\  \ \gm[p+2-k]\ge p\}\Big| =  \big|\hAge\big|.
\end{eqnarray*}

Using these identities and Corollary~\ref{Corollary64|B<|=|A>=|}, we obtain
\begin{align*}
 \big|J_{\Gamma}\big|
   & = \big|\hAbet\big| + 2\big|\hAge\big| + \big|\hBbet\big| + 2\big|\hBge\big|, \\
   & = \big|\hAbet\big| + 2\big|\hAge\big| + \big|\hAbet\big|-2+ 2\left(\big|\hAle\big|-1\right),\\
   & = 2\left(\big|\hAle\big| + \big|\hAbet\big| + \big|\hAge\big|\right)-4,\\
   & = 2\cdot\Big|\Z\cap[2,n+2]\Big|-4 =2n-2.
\end{align*}
\ok

The following corollary comes from proof of the previous theorem.
%%%%%%%%%%%%%%%%%%%%%%%%%%%%%%%%%%%%%%%%%%       Corolario         %%%%%%%%%%%%%%%%
\begin{cor}\label{Corollary65.5}
 Under the hypotheses of Theorem~\ref{Theorem63TheoremAMinusTwo=>Bijective},
\[\Big|\left\{ k\in\Z\ :\ 2\le k\le 2n,\ \ \gk\ge p\right\}\Big|=n,\]
\[\Big|\left\{ k\in\Z\ :\ 2n+1\le k\le 4n,\ \ \gk\ge p\right\}\Big|=n-2.\]
\end{cor}

\section{Sums involving $\r{k^2-k+2-3n}.$}\label{Section:ResiduesInvolvingk^2-k+2-3n}
Consider  $J_n=\Big|\left\{k\in\Z\ :\ 2\le k\le n+2,\ \ \gk\ge p\right\}\Big|$ and we  call $J_n$ simply \textit{the number of jumps}.
We will now develop formulas relating the term residues of $k^2-k+2-3n$ modulus $p.$

%%%%%%%%%%%%%%%%%%%%%%%%%%%%%%%%%%%%%%%%%%       Theorem           %%%%%%%%%%%%%%%%
\begin{thm}\label{Theorem66}
If $n>3$\, and \,$p=4n-1$ is prime then
\begin{eqnarray*}
    \sum_{k=1}^{p-1}\r{k^2-k+2-3n} & = & \sum_{k=1}^{p-1}\r{k^2}+3n-2.\\
    \sum_{k=1}^{2n}\r{k^2-k+2-3n}   & = & \sum_{k=1}^{2n}\r{k^2}+n.\\
    \sum_{k=1}^{n}\r{k^2-k+2-3n}   & = & \sum_{k=1}^{n}\r{k^2}+\frac{n(n+1)}{2}-p(J_n-1-M).
\end{eqnarray*}
\end{thm}
%%%%%%%%%%%%%%%%%%%%%%%%%%%%%%%%%%%%%%%%%%       Proof of Theorem    %%%%%%%%%%%%%%
\noindent {\bf Proof.}
Notice that
\[\defi{\r{x+y}}{\r{x}+\r{y} & \mytext{if} \r{x}+\r{y}<p,\\ \r{x}+\r{y}-p & \mytext{if} \r{x}+\r{y}\ge p.}
\]

Recall that we defined $k$ as a jump when $\gk=\rkm+\r{k+1-3n}\ge p.$ \ By Theorem~\ref{Theorem65},

\begin{eqnarray*}
   \sum_{k=2}^{4n-1} & & \hspace{-.75cm}\r{k^2-k+2-3n} \\
   & = & \sum_{k: \gk\ge p}\hspace{-.2cm}\big(\rkm+\r{k+1-3n}-p\big)\\
   & & \crea{5}{0} +\sum_{k: \gk< p}\hspace{-.2cm}\big(\rkm+\r{k+1-3n}\big),\\
   & = & -p|J_\Gamma|+\sum_{k=2}^{p}\rkm+\sum_{k=2}^{3n-2}\r{k+1-3n}+\sum_{k=3n-1}^{4n-1}\r{k+1-3n},\\
   & = & -p(2n-2)+\sum_{k=1}^{p-1}\rk+\sum_{k=2}^{3n-2}(k+1-3n+p)+\sum_{k=3n-1}^{4n-1}(k+1-3n),\\
   & = &\sum_{k=1}^{p-1}\rk+\sum_{k=2}^{4n-1}(k+1-3n)\ +\ p(2-2n+3n-3), \\
   & = &\sum_{k=1}^{p-1}\rk+ 3n-2.
\end{eqnarray*}
Since $\r{k^2-k+2-3n}=2-3n+p$ for $k=0,p$ the result follows. Observe that
 \[ \sum_{k=0}^{p-1}\r{k^2-k+2-3n} =  \sum_{k=1}^{p-1}\r{k^2}+p.\]
The other two formulas are done analogously.
\ok
Compare this result with Theorem 5.1 in~\cite{garcia:siqrmp22} when $p=4n+1,$
 \[ \sum_{k=0}^{p-1}\r{k^2+k+1-n} =  \sum_{k=1}^{p-1}\r{k^2}-p=p(2n-1).\]

%%%%%%%%%%%%%%%%%%%%%%%%%%%%%%%%%%%%%%%%%%%%%%%%%%%%%%       Remark            %%%%%%%%%%%%%%%%
\begin{rem}\label{Remark44TableSpecialCases}
The first two formulas in Theorem~\ref{Theorem66} are actually valid for any $n$ but the third one is only valid for $n>3.$
 Charts~\ref{TableSpecialExcludedCases} and ~\ref{Table2SpecialExcludedCases} contain  the sum of the residues $\r{k^2-k+2-3n}$ in the special cases excluded in such theorem.
%%%%%%%%%%%%%%%%%%%%%%%%%%%%%%%%%%%%%%%%%%       Table             %%%%%%%%%%%%%%%%
 \begin{table}[h!]
 \centering{\small
\begin{tabular}{|r|c|c|}
  \hline
  % after \\: \hline or \cline{col1-col2} \cline{col3-col4} ...
  $n$ & $\ds{\sum_{k=0}^{p-1}\r{k^2-k+2-3n}}$ &  $\ds{\sum_{k=0}^{p-1}\r{k^2}}+p$  \\ \hline
  $1$ & $5$ & $5$ \\ \hline
  $2$ & $21$ & $21$  \\ \hline
  $3$ & $55$ & $55$  \\ \hline
\end{tabular}
 \caption{Sum of residues $\r{k^2-k+2-3n}$ in special cases.}\label{TableSpecialExcludedCases}}
\end{table}
%%%%%%%%%%%%%%%%
%%%%%%%%%%%%%%%%%%%%%%%%%%%%%%%%%%%%%%%%%%       Table

%%%%%%%%%%%%%%%%%%%%%%%%%%%%%%%%%%%%%%%%%%       Table             %%%%%%%%%%%%%%%%
 \begin{table}[h!]
 \centering{\small
\begin{tabular}{|r|c|c|}
  \hline
  % after \\: \hline or \cline{col1-col2} \cline{col3-col4} ...
  $n$ &  $\ds{\sum_{k=0}^{2n}\r{k^2-k+2-3n}}$ &  $\ds{\sum_{k=0}^{2n}\r{k^2}}+2n+1$  \\ \hline
  $1$ & $5$& $5$ \\ \hline
  $2$ & $14$& $14$ \\ \hline
  $3$ & $32$& $32$  \\ \hline
\end{tabular}
 \caption{Sum of residues $\r{k^2-k+2-3n}$ in special cases.}\label{Table2SpecialExcludedCases}}
\end{table}
%%%%%%%%%%%%%%%%
%%%%%%%%%%%%%%%%%%%%%%%%%%%%%%%%%%%%%%%%%%       Table

\end{rem}

%%%%%%%%%%%%%%%%%%%%%%%%%%%%%%%%%%%%%%%%%%       Corolario         %%%%%%%%%%%%%%%%
\begin{cor}\label{Theorem67}
If $n>3$\, and \,$p=4n-1$ is prime then

{\small
\begin{eqnarray*}
    \sum_{k=1}^{p-1}\ffloor{k^2-k+2-3n}{p} & = & \frac{p(p-5)+6-n}{3}-\frac{1}{p}\sum_{k=1}^{p-1}\rk.\\
   \sum_{k=1}^{2n}\ffloor{k^2-k+2-3n}{p} & = & \frac{n(2n-1)(4n-7)}{3p}-\frac{1}{p}\sum_{k=1}^{2n}\rk.\\
   \sum_{k=1}^{n}\ffloor{k^2-k+2-3n}{p} & = & \frac{(n+1)(2n^2-23n+6)}{3p}-\frac{1}{p}\sum_{k=1}^{n}\rk+J_n-M.
\end{eqnarray*} }
\end{cor}
%%%%%%%%%%%%%%%%%%%%%%%%%%%%%%%%%%%%%%%%%%       Proof of Corollary  %%%%%%%%%%%%%%
\noindent {\bf Proof.}
From $x=p\ffloor{x}{p}+\r{x},$ we obtain that if $y=\sum_{k=1}^{p-1}\ffloor{k^2-k+2-3n}{p}$ then
\[\sum_{k=1}^{p-1} \big(k^2-k+2-3n\big) =py+\sum_{k=1}^{p-1}\r{k^2-k+2-3n}.\]
From Theorem~\ref{Theorem66},
\[\frac{(p-1)p(2p-1)}{6}-\frac{(p-1)p}{2}+(p-1)(2-3n)=py+\sum_{k=1}^{p-1}\rk+3n-2,\]
then
\[p\cdot\frac{p\cdot(p-5)+6-n}{3}=py+\sum_{k=1}^{p-1}\rk.\]
The result follows. The proofs of the other two sums are done similarly.
\ok

The purpose of the following lemma and remark is to find a formula for $\sum_{k=0}^{n}\ffloor{k^2-k+2-3n}{p}.$

%%%%%%%%%%%%%%%%%%%%%%%%%%%%%%%%%%%%%%%%%%       Lemma             %%%%%%%%%%%%%%%%
\begin{lem}\label{Lemma68TrueMeaningofJumps}
Let $\ds{Q_m=\yp{m}}$ as in Notation~\ref{NotationQmRmM}. Then
\[ \floor{Q_m}=\max\{k\in\N\ :\ \, k^2-k+2-3n\le mp\}. \]
\end{lem}
%%%%%%%%%%%%%%%%%%%%%%%%%%%%%%%%%%%%%%%%%%       Proof of Lemma    %%%%%%%%%%%%%%%%
\noindent {\bf Proof.}

Clearly $\floor{Q_m}\le Q_m<\floor{Q_m}+1$ and since $Q_m$ is the non-negative root of $x^2-x+2-3n=mp, \ k_0=\floor{Q_m}$ satisfies $0\le k_0\le Q_m$ and $k_0^2-k_0+2-3n\le mp.$ Hence $k_1=\floor{Q_m}+1$ satisfies $0\le k_1 $ and $k_1^2-k_1+2-3n>mp.$ If $k_0=0$ then $k_1=1$ and hence $p\le mp<2-3n$ which is impossible, therefore $k_0\in\N.$
\ok

%%%%%%%%%%%%%%%%%%%%%%%%%%%%%%%%%%%%%%%%%%       Remark            %%%%%%%%%%%%%%%%
\begin{obs}\label{Observation69}
By Lemma~\ref{Lemma53QhInequality}, there are no integers $k, m, n$ such that $\floor{Q_m}=Q_m,$ i.e. $k^2-k+2-3n\neq mp$ regardless of $k,m,n.$
\end{obs}

%%%%%%%%%%%%%%%%%%%%%%%%%%%%%%%%%%%%%%%%%%       Theorem           %%%%%%%%%%%%%%%%
\begin{thm}\label{Theorem70}
  Let $n> 3$ and $p=4n-1$ prime.  Then
  \[ \sum_{k=1}^{n}\ffloor{k^2-k+2-3n}{p}=(M-1)n-\sum_{m=0}^{M-1}\floor{Q_m}.\]
\end{thm}
%%%%%%%%%%%%%%%%%%%%%%%%%%%%%%%%%%%%%%%%%%       Proof of Theorem    %%%%%%%%%%%%%%
\noindent {\bf Proof.}
Since $n>3, 1\le M.$ Define, $t_m=\floor{Q_m}$ and
\[ \defi{H_m}{\big\{\,1,...,t_0\big\} & \mytext{if} m=0,  \\[1mm]
              \big\{t_{m-1}+1,t_{m-1}+2,...,t_m\big\}
              & \mytext{if} 1\le m\le M-1,  \\[1mm]
              \big\{t_{M-1}+1,t_{M-1}+2,...,n\big\}
              &\mytext{if} m=M.}
\]

By Lemma~\ref{Lemma68TrueMeaningofJumps}, $k=\floor{Q_m}$ satisfies $k^2-k+2-3n\le mp$ and from Observation~\ref{Observation69}, $k^2-k+2-3n< mp.$ Therefore  by Lemma~\ref{Lemma68TrueMeaningofJumps}, for $k\in H_m$ we have
\[(m-1)p<k^2-k+2-3n<mp.\]

For such $k$ necessarily $\floor{(k^2-k+2-3n)/p}=m-1.$

Notice that $\{H_0,H_1,...,H_M\}$ is a partition of $\{1,2...,n\}$ as $Q_{M-1}<n<Q_M$ (see Lemma 2.5 in~\cite{garcia:cfcniqf21}).   Therefore
        \begin{eqnarray*}
          \sum_{k=1}^{n}\ffloor{k^2-k+2-3n}{p} &=& \sum_{m=0}^{M} \sum_{\substack{k=0\\ k\in H_m}}^{n}\ffloor{k^2-k+2-3n}{p}, \\
           &=& \sum_{m=0}^{M} (m-1)|H_m|, \\
           &=& -t_0+0 (t_1-t_0)+
            1 (t_2-t_1)+\cdots,\\
            & &+(M-2) (t_{M-1}-t_{M-2})
            +(M-1) (n-t_{M-1}),\\
            &=& -\sum_{m=0}^{M-1}\floor{Q_m}+(M-1)n.
        \end{eqnarray*}\ok

%%%%%%%%%%%%%%%%%%%%%%%%%%%%%%%%%%%%%%%%%%       Corolario         %%%%%%%%%%%%%%%%
\begin{cor}\label{Corollary71}
If $n>3$\, and \,$p=4n-1$ is prime then

\begin{equation*}
   \sum_{k=0}^{n}\ffloor{k^2}{p}+\sum_{k=1}^{n}\ffloor{k^2-k+2-3n}{p}=\frac{n(n-5)}{6}+M-\frac{1}{2p}\sum_{k=1}^{p-1}\r{k^2}.
\end{equation*}

\end{cor}
%%%%%%%%%%%%%%%%%%%%%%%%%%%%%%%%%%%%%%%%%%       Proof of Corollary  %%%%%%%%%%%%%%
\noindent {\bf Proof.}
From~\cite{zeller:usvggb79} (page 253) we have
\begin{equation}\label{71-1}
\sum_{m=1}^{M}\floor{R_m}=Mn-\sum_{k=0}^{n}\ffloor{k^2}{p}.
\end{equation}

\nn Corollary 2 in~\cite{garcia:cfcniqf21} states
\begin{equation}\label{71-2}
     \frac12\sum_{k=1}^{p-1}\r{k^2} = p\left(\sum_{m=1}^{M} \frm +\sum_{m=0}^{M-1} \fqm\right)-Mp(2n-1)+\frac{p\cdot(n^2+n)}{6}.
  \end{equation}
\nn Using Theorem~\ref{Theorem70} and Equation~\ref{71-1} in Equation~\ref{71-2} we obtain
        \begin{eqnarray*}
          \frac12\sum_{k=1}^{p-1}\r{k^2} & = & p\left((2M-1)n-\sum_{k=1}^{n}\ffloor{k^2}{p} -\sum_{k=1}^{n}\ffloor{k^2-k+2-3n}{p} \right)\\
          & & -Mp(2n-1)+\frac{p\cdot(n^2+n)}{6}.
        \end{eqnarray*}
The result now follows.\ok

Compare Corollary~\ref{Corollary71} with Theorem 2.2 (case $p=4n+1$) in~\cite{garcia:siqrmp22} which can be rewritten as
\begin{equation*}
  \sum_{k=0}^{n}\ffloor{k^2}{p}+\sum_{k=1}^{n}\ffloor{k^2+k+1-n}{p}=\frac{(n+3)(n+2)}{6}+M-\frac{1}{2p}\sum_{k=1}^{p-1}\r{k^2}+u_n.
\end{equation*}

\bigskip

%%%%%%%%%%%%%%%%%%%%%%%%%%%%%%%%%%%%%%%%%%       Theorem           %%%%%%%%%%%%%%%%
\begin{cor}\label{Theorem73}
If $n>3$\, and \,$p=4n-1$ is prime then
  \[ \sum_{m=1}^{M} \fr{m}- \sum_{m=0}^{M-1}\fq{m}= J_n-M-1.\]
\end{cor}

%%%%%%%%%%%%%%%%%%%%%%%%%%%%%%%%%%%%%%%%%%       Proof of Theorem    %%%%%%%%%%%%%%
\noindent {\bf Proof.}
From Lemmas~\ref{Lemma57ThreeChunks},~\ref{Lemma62Bbet} and~\ref{Lemma72} we obtain

\begin{eqnarray*}
  J_n &=& |\hAbet|+|\hAge| = |\hBbet|+2+|\hAge|, \\
    &=& M_0+2 +\sum_{m=0}^{M_0}\left(\ell_m- k_m\right),\\
    &=& M+1 +\sum_{m=0}^{M_0}\left(\floor{R_{m+1}}- \floor{Q_m}\right),\\
    &=& M+1 +\sum_{m=1}^{M} \fr{m}- \sum_{m=0}^{M-1}\fq{m}.
\end{eqnarray*}\ok

Compare Corollary~\ref{Theorem73} with Theorem 5.3 (case $p=4n+1$) in~\cite{garcia:siqrmp22}
\begin{equation*}
  \sum_{m=1}^{M}\fr{m}- \sum_{m=0}^{M}\fs{m}  =j_n+2-n-u_n.
\end{equation*}

\section{Class Number Identities}\label{Section:ClassNumber}
   In this section, we establish some identities involving the class number $h=h(-p)$ of the imaginary quadratic field $\Q(\sqrt{-p})$ when $p$ is of the form $p=4n-1.$ These identities are based on the previous formulas we have developed in previous sections.

   In~\cite{garcia:cfcniqf21}, we have

\begin{eqnarray*}
    h &=& (2M+1)(2n-1)-2\left(\sum_{m=1}^{M} \fr{m}- \sum_{m=0}^{M-1}\fq{m}\right)-\frac{n^2+n}{3}.\\
    h &=& \frac{p-1}{2}-\frac{1}{p}\sum_{k=1}^{p-1}\rk.
\end{eqnarray*}

\bigskip

\nn From Corollaries~\ref{Theorem67},~\ref{Corollary71} and $\sum_{k=1}^{p-1}\rk=2\sum_{k=1}^{2n}\rk-2n$ we obtain
\begin{cor}\label{Corollary74}
\begin{eqnarray*}
   %\sum_{k=0}^{n}\ffloor{k^2-k+2-3n}{p} &=& \frac{h}{2}+\sum_{k=0}^{M}\fr{k}+M(1-n)+\frac{n^2-11n-3}{6},\\
   \sum_{k=1}^{p-1}\ffloor{k^2-k+2-3n}{p} &=& h+\frac{16n^2-35n+15}{3}.\\
   \sum_{k=1}^{2n}\ffloor{k^2-k+2-3n}{p} &=& \frac{h}{2}+\frac{4n^2-14n+9}{6}.\\
   \sum_{k=1}^{n}\ffloor{k^2-k+2-3n}{p} &=& \frac{h}{4}+\frac{J_n}{2}+\frac{n^2-17n-3}{12}.\\
  \sum_{k=1}^{n}\ffloor{k^2-k+2-3n}{p} &=& \frac{h}{2}+M+\frac{n^2-11n+3}{6}- \sum_{k=1}^{n}\ffloor{k^2}{p}.\\
  \sum_{k=1}^{n}\ffloor{k^2-k+2-3n}{p} &=& \frac{h}{2}+\sum_{m=1}^{M}\frm+M(1-n)+\frac{n^2-11n+3}{6}.
 \end{eqnarray*}
\end{cor}

\bigskip

\nn  From Corollary~\ref{Theorem73}, we have
\begin{cor}\label{Corollary75}
\begin{eqnarray*}
   \frac{J_n}{2}+M(n-1) & = &  \frac{h}{4}+\sum_{m=1}^{M}\frm+\frac{n^2-5n+9}{12}. \\
   -\frac{J_n}{2}+Mn & = &  \frac{h}{4}+\sum_{m=0}^{M-1}\fqm+\frac{n^2-5n-3}{12}.\\
   \frac{J_n}{2}+\sum_{k=1}^{n-1}\ffloor{k^2}{p} & = &  \frac{h}{4}+\frac{n^2-5n+9}{12}.
 \end{eqnarray*}

\end{cor}

\bigskip

\nn   From Theorem~\ref{Theorem66}, we conclude
\begin{cor}\label{Corollary76}
\begin{eqnarray*}
    \sum_{k=1}^{p-1}\r{k^2-k+2-3n} & = & p\cdot\left(2n-h\right)-n-1.\\
    \sum_{k=1}^{2n}\r{k^2-k+2-3n}   & = & \frac{p\cdot(2n-h)+1}{2}.\\
    \sum_{k=1}^{n}\r{k^2-k+2-3n}   & = & \frac{p}{4}(3n+2-2J_n-h)-\frac{(n+1)(n-1)}{4}.
\end{eqnarray*}
\end{cor}

\bigskip

\nn Finally, combining this last formula with Theorem~\ref{Theorem66} we have
\begin{cor}\label{Corollary77}
\[\sum_{k=1}^{n}\r{k^2}  = \frac{p}{4}(2J_n+2n-3-4M-h)+\frac{n(n+1)}{4}.\]
\end{cor}

\bigskip

\nn The numerical data we have allow us to pose the following
%%%%%%%%%%%%%%%%%%%%%%%%%%%%%%%%%%%%%%%%%%       Conjecture           %%%%%%%%%%%%%%%%
\begin{cnjt}\label{conjecture78}
  Consider all $n$ such that $p=4n-1$ a prime number.
  Then
\[\lim_{n\to\infty} \frac{J_n}{n}= \frac{3}{8}.\]
\end{cnjt}

\medskip

%%%%%%%%%%%%%%%%%%%%%%%%%%%%%%%%%%%%%%%%%%       Conjecture           %%%%%%%%%%%%%%%%
\begin{cnjt}\label{conjecture79}
  Consider all $n$ such that $p=4n-1$ a prime number.
  Then
\[\lim_{n\to\infty} \frac{\sum_{m=1}^{M}\frm+ \sum_{m=0}^{M-1}\fqm}{Mp+2n}  = \frac{1}{3}.\]
 Also a good estimate of $\sum_{m=1}^{M}\frm+ \sum_{m=0}^{M-1}\fqm$ is $\floor{(Mp+2n)/3}.$
\end{cnjt}

\medskip

%For instance when $n=200,000,\ \sum_{k=1}^{M}\floor{R_{k}}+ \sum_{k=0}^{M-1}\floor{Q_k}=13,333,449,738$ and $\floor{(Mp+2n)/3}= 13,333,450,000,$ in fact the error is  $262.$

%%%%%%%%%%%%%%%%%%%%%%%%%%%%%%%%%%%%%%%%%%       Conjecture           %%%%%%%%%%%%%%%%
%\begin{cnjt}\label{conjecture80}
%  Consider all $n$ such that $p=4n-1$ a prime number.
%  Then
%\[\lim_{n\to\infty} \frac{h(-p)}{n^2}  = 0.\]
%\end{cnjt}

\newpage
 \bibliographystyle{siam}
 \bibliography{quady}
\end{document}